\documentclass[a4paper,twoside]{article}

\usepackage{graphicx}
\usepackage{epsfig}
\usepackage{subcaption}
\usepackage{calc}
\usepackage{amssymb}
\usepackage{amstext}
\usepackage{amsmath}
\usepackage{amsthm}
\usepackage{multicol}
\usepackage{pslatex}
\usepackage{apalike}
\usepackage{here}
\usepackage[ruled,vlined]{algorithm2e}
\usepackage[bottom]{footmisc}
\usepackage{SCITEPRESS}     
\theoremstyle{definition}

\newtheorem{Theorem}{Theorem}[section]

\newtheorem{Proof}{Proof}

\renewcommand{\abstractname}{Abstract}
\begin{document}
\title{Solving Monge problem by \\Hilbert space embeddings of probability measures}

\author{\authorname{Takafumi Saito\sup{1}
and Yumiharu Nakano\sup{1}
}
\affiliation{\sup{1}Department of Mathematical and Computing Science, Institute of Science Tokyo, Tokyo, Japan}
\email{saito.t.8e06@m.isct.ac.jp, nakano.y.ee47@m.isct.ac.jp}
}


\abstract{We propose deep learning methods for classical Monge's optimal mass transportation problems, where where the distribution constraint is treated as a penalty term defined by the maximum mean discrepancy in the theory of Hilbert space embeddings of probability measures. We prove that the transport maps given by the proposed methods converge to optimal transport maps in the problem with $L^2$ cost. Several numerical experiments validate our methods. In particular, we show that our methods are applicable to large-scale Monge problems. This is a corrected version of the ICORES 2025 proceedings paper.}

\onecolumn \maketitle \normalsize \setcounter{footnote}{0} \vfill

\section{\uppercase{Introduction}}
\label{sec:introduction}
Our aim in this paper is to propose numerical methods for Monge's mass transportation problem, described as follows: given two Borel
probability measures $\mu$, $\nu$ on $\mathbb{R}^d$ and a cost function $c:\mathbb{R}^d\times\mathbb{R}^d\to [0,\infty]$, to minimize
Our problem is to minimize 
\[
 M(T):=\int_{\mathbb{R}^d}c(x,T(x))d\mu(x)
\]
over all Borel measurable mapping $T:\mathbb{R}^d\to\mathbb{R}^d$ such that 
$\mu\circ T^{-1}=\nu$. Here, $\mu\circ T^{-1}$ denotes the pushforward of $\mu$ with $T$, i.e., $\mu\circ T^{-1}(A)=\mu(T\in A)$ for any Borel set $A$. 

We shall briefly describe the background of the Monge problem. Monge problem was proposed by Gaspard Monge in 1781 \cite{M}. 
In the 20th century, this problem was expanded by Kantorovich to make it mathematically easier to handle, and is now called the Monge-Kantrovich problem \cite{K1}, \cite{K2}. 

The Monge problem, as extended by Kantrovich, is called the Monge-Kantrovich problem. 
The first existence and uniqueness result was established by Brein, \cite{b1}, \cite{b2}. 
Gangbo and McCann further developed the generalized problem \cite{G}.
Mikami also provided a probabilistic proof of Breiner's result \cite{M2}.

Algorithms for solving the Monge-Kantrovich problem can be traced back nearly 100 years \cite{To}. 
Since the advent of mathematical programming, these algorithms have been a field of significant interest \cite{Da2}. 
This is largely because Dantzig's initial motivation for developing mathematical programming was related to solving transportation problems \cite{Da1}, 
and it was later discovered that optimal transportation problems and minimum cost flow problems are equivalent \cite{Co}. 
Research has advanced in the development of solvers that apply mathematical programming to machine learning \cite{Bau}, 
as well as in solvers used in dynamic optimal transport (OT) \cite{Papa}. 
Today, it remains one of the most actively researched fields.

We describe recent solvers developed to address the Monge-Kantrovich problem. In recent years, the Monge-Kantrovich problem has been actively applied in machine learning. One prominent approach involves solving the entropy optimization problem by introducing an entropy regularization term, using the Sinkhorn algorithm \cite{cu}. The Python library {\tt POT} \cite{P} provides a solver for various optimal transport problems, including entropy-regularized optimal transport.

When solving the optimal transport problem, computations can scale with the cube of the input data size. Therefore, it is crucial to support large-scale computations and to leverage GPUs for numerical calculations. Optimal transport solvers designed for large-scale computations include the Python library {\tt Geomloss} \cite{Geom}. Another tool, the Python library {\tt KeOps} \cite{keO}, optimizes the reduction of large arrays using neural networks and kernel formulas.

Additionally, the {\tt OTT} library \cite{OTT}, based on the high-performance JAX library for numerical calculations, offers implementations that can solve problems such as the Monge-Kantrovich problem.

Our aim is to derive a numerical solution for the basic mathematical analysis of the original Monge problem, rather than the Monge-Kantrovich problem. In particular, the goal is to develop algorithns capable of performing GPU-based numerical computations and handling large-scale calculations. The method for deriving the numerical solution primarily utilizes the embedding theory of probability measures, which was previously applied by Nakano \cite{N} to obtain a numerical solution for the Schr{\"o}dinger bridge problem. In this study, we also apply this theory to derive a numerical solution for the Monge problem. The penalty method is employed to find the optimal solution, with the use of Maximum Mean Discrepancy (MMD) as the penalty function being a novel approach. Unlike existing methods, this method is independent of input data size. 
We also confirm this through numerical experiments. 

This paper is structured as follows. In the next section, we review some basic results on the theory of Hilbert space embedding of probability measures and describe a numerical method that describes a theoretical approximation solution to the Monge problem with $L^2$ cost. 
Section \ref{sec:3} gives numerical experiments.

\section{\uppercase{PENALTY METHOD}}\label{sec:2}
\subsection{Hilbert space embeddings of probability measures}\label{sec:2.1}
We shall give a quick review of theory of Hilbert space embeddings of probability measures, as developed in Sriperumbudur et al. \cite{S}. 
Denote by $\mathcal{P}(\mathbb{R}^d)$ the set of all Borel probability measures on $\mathbb{R}^d$. 
Let $K$ be a symmetric and 
 strictly positive definite kernel on $\mathbb{R}^{d}$, i.e., 
 $K(x,y)=K(y,x)$ for 
 $x,y\in\mathbb{R}^{d}$ and 
  for any positive distinct $x_{1},\ldots,x_{N}\in\mathbb{R}^{d}$ and 
  $\alpha=(\alpha_{1},\ldots,\alpha_{N})  ^{\mathrm{T}}\in\mathbb{R}
  ^{N}\backslash\{0\}$,
\begin{equation*}
\sum^{N}_{j,\ell=1}
\alpha_{j}\alpha_{\ell}K(x_{j},x_{\ell})>0.
\end{equation*}
Assume further that $K$ is bounded and continuous on $\mathbb{R}^d\times\mathbb{R}^d$. 
Then, there exists a unique Hilbert space 
$\mathcal{H}$
 such that $K$ is a reproducing kernel on 
 $\mathcal{H}$ with norm $\|\cdot\|$ \cite{wen:2010}. 
 Then consider the {\it maximum mean discrepancy (MMD)} $\gamma_{K}$ defined by, 
 for 
 $\mu_{0},\mu_{1}\in\mathcal{P}(\mathbb{R}^{d})$,
\begin{align*}
 \gamma_{K}(\mu_{0},\mu_{1}):
 &=
 \mathrm{sup}_{f\in\mathcal{H},\|f\|\leqq 1}
 \left|\int_{\mathbb{R}^{d}}f d\mu_{0}- \int_{\mathbb{R}^{d}}fd\mu_{1}\right|\notag\\
 &=\left\|\int_{\mathbb{R}^{d}}K(\cdot,x)\mu_{0}(dx)-
 \int_{\mathbb{R}^{d}}K(\cdot,x)\mu_{1}(dx)\right\|_{\mathcal{H}} 
 \end{align*}
It is known that if $K$ is an integrally strictly positive definite then $\gamma_K$ defines a metric on $\mathcal{P}(\mathbb{R}^d)$. 
Examples of integrally strictly positive definite kernels include the Gaussian kernel 
$K(x,y)=e^{-\alpha |x-y|^2}$, $x,y\in\mathbb{R}^d$, where $\alpha>0$ is a constant, and 
the Mat{\'e}rn kernel $K(x,y)=K_{\alpha}(|x-y|)$, $x,y\in\mathbb{R}^d$, where $K_{\alpha}$ is the modified Bessel function of 
order $\alpha>0$. It is also known that Gaussian kernel as well as Mat{\'e}rn kernel metrize the weak topology on $\mathcal{P}(\mathbb{R}^d)$ 
\cite{S}, \cite{N}. 

Define
\begin{align*}
K_{1}(x,y):&=K(x,y)-\int_{\mathbb{R}^{d}}K(x,y^{\prime})
\mu_{1}(dy^{\prime})\\
&\quad -\int_{\mathbb{R}^{d}} K(x^{\prime},y)\mu_{1}(dx^{\prime}) .
\end{align*}
Then, 
 \begin{equation}
 \label{eq:1}
 \begin{aligned}
 \gamma_{K}(\mu_{0},\mu_{1})^{2}&=
 \int_{\mathbb{R}^{d}}\int_{\mathbb{R}^{d}}
 K_{1}(x,y)\mu_{0}(dx)\mu_{0}(dy)\\ 
 &\quad +
 \int_{\mathbb{R}^{d}\times\mathbb{R}^{d}}
 K(x^{\prime},y^{\prime})\mu_{1}(dx^{\prime})\mu_{1}(dy^{\prime}).
 \end{aligned}
 \end{equation}
\subsection{Proposed methods}
Let $\gamma=\gamma_K$ be as in Section \ref{sec:2.1}. 
We proposed a penalty method for Monge problem by Hilbert space embeddings of probability
 measures, described as follows:
\begin{equation*}
\inf_{T} M_{\lambda}(T):=\inf_{T}\left\{M(T) + \lambda \gamma^{2}(\mu\circ T^{-1}, \nu)\right\}.
\end{equation*}
Note that the second term penalizes the distance between the laws of $T(x)$ and $\nu$.
Moreover,
the second term in the above formula can be expressed discretely as follows: given IID samples
$X_{1},\cdots,X_{M}\sim\mu_{0}$ and 
$Y_{1},\cdots,Y_{M}\sim\mu_{1}$, an unbiased estimator of 
$\gamma_{k}(\mu_{0},\mu_{1})$ is given by
\begin{equation}
\label{eq:2}
\begin{aligned}
\bar{\gamma}_{K}(\mu_{0},\mu_{1})^{2}&:=
\frac{1}{M(M-1)}\sum_{i}\sum_{j\neq i}K(X_{i},X_{j})
\\
&\quad -\frac{2}{M^{2}}\sum_{i,j}K(X_{i},Y_{j})
\\
&\quad +\frac{1}{M(M-1)}\sum_{i}\sum_{j \neq i}K(Y_{i},Y_{j}) 
\end{aligned}
\end{equation}
\cite{GBRSS}. 
Then, we approxiate $T$ by a class 
$\{T_{\theta}\}_{\theta\in\Theta}$ of deep neural networks. 
Each $T_{\theta}$ can be given by a multilayer perception with input layer $g_{0}$, 
 $L-1$ hidden layer $g_{1},\ldots,g_{L-1}$, and output layer $g_{L}$, where $L\geqq1$ 
 and for $\xi\in\mathbb{R}^{1+m}$, $g_{0}(\xi)=\xi$, $g_{\ell}(\xi)=\phi_{\ell-1}(w_{\ell}g_{\ell-1}(\xi)+\beta_{\ell}) \in\mathbb{R}^{m_{\ell}}, \ell=1,\ldots,L$,  
 for some matrices $w_{\ell}$ and vectors $\beta_{\ell},\ell=1,\ldots,L$.
  Here $m_{\ell}$ denotes the number of units in the layer $\ell$, and $\phi_{\ell-}$ is
   an activation function. Then the parameter $\theta$ is described by $\theta=(w_{\ell},
   \beta_{\ell})_{\ell=1,\ldots,L}$ 
   and $T_{\theta}(x)=g_{L}(x)\in\mathbb{R}^{d}$. 
   For $\lambda>0$, the integral term in \eqref{eq:1} is replaced by \eqref{eq:2} and 
   $T_{\theta}(x)$ as follows: by Subsetcion 3.2. in \cite{N},
\begin{equation}
\label{eq:3}
\begin{aligned}
F_{1}(\theta)&=\frac{1}{\lambda M}\sum_{i=1}^MC(X_i,T_{\theta}(X_i)) \\ 
&\quad +\frac{1}{M(M-1)}\sum_{i}\sum_{j\neq i}K(T_{\theta}(X_{i}),
T_{\theta}(X_{j}))\\
&\quad -\frac{2}{M^{2}}K(T_{\theta}(X_{i}),Y_{j})
\end{aligned}
\end{equation}
The algorithm is described below and we test our one algorithm thorough a numerical experiment. 
\begin{algorithm}[h]
 \caption{Deep learning algorithm with empirical MMD}
 \DontPrintSemicolon
  \KwData{The number $n$ of the iterations, the batch size $M$, weight parameter $\lambda>0$}
 \KwResult{transport map $T$ }
 initialization\;
$X_{1},\ldots,X_{M}
    \leftarrow~\mathrm{IID}$ samples from  $\mu$,\;
    $Y_{1},\cdots,Y_{M}
    \leftarrow~\mathrm{IID}$ samples from $\nu$.\;
    \For {k=1,2,$\ldots,n$}{
Compute $F_{1}(\theta)$ in $(3)$ using
     $\{X_{j}, Y_{j}\}_{j=1}^M$. \; 
Take the gradient step on $\nabla_{\theta}F_{1}(\theta)$.
 }
\end{algorithm}
\subsection{Theoretical result}

For given $\lambda>0$ consider the minimization problem  
\begin{equation*} 
 M_{\lambda}(T):= M(T) + \lambda \gamma^{2}(\mu\circ T^{-1}, \nu)
\end{equation*}
over all Borel measurable mappings $T$. 
Take arbitrary positive sequences $\{\varepsilon_n\}_{n=1}^{\infty}$ and 
$\{\lambda_n\}_{n=1}^{\infty}$ such that 
\begin{equation*}
 \varepsilon_n\to 0, \quad \lambda_n\to +\infty \quad (n\to \infty). 
\end{equation*}
Then take $T_n:\mathbb{R}^d\to\mathbb{R}^d$ such that 
\begin{equation*}
 M_{\lambda_n}(T_n)-\varepsilon_n\le M^*_{\lambda_n}:=\inf_TM_{\lambda_n}(T).  
\end{equation*}
Then we have the following:
\begin{Theorem}
Let $c(x,y)=|x-y|^2$. Suppose that $\mu$ is absolutely continuous with respect to the Lebesgue measure and that 
\begin{equation*}
 \int_{\mathbb{R}^d}|x|^2\mu(dx) + \int_{\mathbb{R}^d}|y|^2\nu(dy)< \infty. 
\end{equation*}
Suppose moreover that $\gamma$ metrizes the weak topology on $\mathcal{P}(\mathbb{R}^d)$. Then, 
\begin{align}
 &\lim_{n\to\infty}\sqrt{\lambda_n}\gamma(\mu\circ T_n^{-1},\nu) =0, \label{eq:4}
 \\ 
 &\lim_{n\to\infty}M_{\lambda_n}(T_n)=M(T^*),  \label{eq:5}
\end{align}
where $T^*$ is the unique optimal transport map. 
In particular, $\{T_n\}_{n=1}^{\infty}$ converges to $T^*$ in law under $\mu$. 
\end{Theorem}
\begin{Proof}
First, note that under our assumption, an optimal transport map does exist uniquely (see, e.g., Theorem 2.12 in \cite{Vi} and Theorem 1.25 in \cite{Sa}).

Step (i). We will show \eqref{eq:4}. 
This claim can be proved by almost the same argument as that given in the proof of Theorem 3.1 in \cite{N}, but we shall give a proof for reader's convenience. 
Assume contrary that 
\begin{equation*} 
 \limsup_{n\to\infty}\lambda_n\gamma^{2}(\mu\circ T_n^{-1},\nu)=5\delta
\end{equation*}
for some $\delta>0$. Then there exists a subsequence $\{n_k\}$ such that 
\begin{equation*}
 \lim_{k\to\infty}\lambda_{n_k}\gamma^{2}(\mu\circ T_{n_k}^{-1},\nu)=5\delta. 
\end{equation*}
Since $\gamma$ is a metric, we have $\gamma(\mu\circ (T^*)^{-1},\nu)=0$, 
whence $M^*_{\lambda}\le M(T^*)$ for any $\lambda$. This means 
\begin{equation*}
 M_{\lambda_n}(T_n)\le M^*_{\lambda_n}+\varepsilon_n \le M(T^*)+\varepsilon_n. 
\end{equation*}
Thus, the sequence $\{M_{\lambda_n}(T_n)\}_{n=1}^{\infty}$ is bounded, and so 
there exists a further subsequence $\{n_{k_m}\}$ such that 
\begin{equation*}
 \lim_{m\to\infty}M_{\bar{\lambda}_m}(\bar{T}_m)=\kappa:=\limsup_{k\to\infty} M_{\lambda_{n_k}}(T_{n_k}) < \infty, 
\end{equation*}
where $\bar{\lambda}_m=\lambda_{n_{k_m}}$ and $\bar{T}_m=T_{n_{k_m}}$. 
Now choose $m_0$ and $m_1$ such that 
\begin{gather*}
\kappa< M_{\bar{\lambda}_{m_0}}(\bar{T}_{m_0}) + \delta,\;\; 
  M_{\bar{\lambda}_{m_1}}(\bar{T}_{m_1}) < \kappa +  \delta, \\   
 \bar{\lambda}_{m_1}> 7\bar{\lambda}_{m_0}, \;\; 
 3\delta + \bar{\varepsilon}_{m_0}< \bar{\lambda}_{m_1}\bar{\gamma}_{m_1}^2< 7\delta, 
\end{gather*}
where $\bar{\gamma}_m=\gamma(\mu\circ T_{n_{k_m}},\nu)$. 
With these choices it follows that 
\begin{align*}
 \kappa &< M_{\bar{\lambda}_{m_0}}(\bar{T}_{m_0}) + \delta 
 \le M^*_{\bar{\lambda}_{m_0}} + \bar{\varepsilon}_{m_0} + \delta \\ 
 &\le M_{\bar{\lambda}_{m_0}}(\bar{T}_{m_1}) + \bar{\varepsilon}_{m_0} + \delta \notag\\ 
 &= M(\bar{T}_{m_1}) + \frac{\bar{\lambda}_{m_0}}{\bar{\lambda}_{m_1}}\bar{\lambda}_{m_1}\bar{\gamma}_{m_1}^2 + \bar{\varepsilon}_{m_0} + \delta \notag\\ 
 &< M(\bar{T}_{m_1}) + \frac{1}{7}\bar{\lambda}_{m_1}\bar{\gamma}_{m_1}^2 + \bar{\varepsilon}_{m_0} + \delta \\
 &< M(\bar{T}_{m_1}) + 2\delta + \bar{\varepsilon}_{m_0} \notag\\ 
 &< M(\bar{T}_{m_1}) + \bar{\lambda}_{m_1}\bar{\gamma}_{m_1}^2 - \delta 
 = M_{\bar{\lambda}_{m_1}}(\bar{T}_{m_1}) - \delta < \kappa, 
\end{align*}
which is impossible.  

Step (ii). Next we will show \eqref{eq:5}. Let $\varepsilon^{\prime}>0$. 
For $R>0$ we have 
\begin{align*}
 \mu\circ T_n(|x|>R)\le&\frac{1}{R^2}\int_{\mathbb{R}^d}|T_n(x)|^2 \mu(dx) 
 \notag\\
 \le&\frac{2}{R^2}\int_{\mathbb{R}^d}|x-T_n(x)|^2\mu(dx) \notag\\ 
  &+ \frac{2}{R^2}\int_{\mathbb{R}^d}|x|^2\mu(dx) \notag\\ 
  &\le \frac{2}{R^2}(M(T^{\ast})+\varepsilon_{n}) + \frac{C}{R^2}
\end{align*}
for some constant $C>0$. Thus we can take a sufficiently large $R^{\prime}$ such that $\mu_n\circ T_n(|x|>R^{\prime})\le \varepsilon^{\prime}$. 
This means that $\{\mu\circ T_n^{-1}\}_{n=1}^{\infty}$ is tight, whence 
there exists a subsequence $\{n_k\}$ such that 
$\mu\circ T_{n_k}^{-1}$ converges weakly to some $\mu^*$. Since we have assumed that $\gamma$ metrizes the weak topology, we get 
\begin{equation*}
 \lim_{k\to\infty}\gamma(\mu\circ T_{n_k}^{-1},\mu^*)=0.
\end{equation*}
This together with the step (i) yields 
\begin{equation*}
 \gamma(\mu^*,\nu)\le \gamma(\mu^*,\mu\circ T_{n_k}^{-1}) 
 + \gamma(\mu\circ T_{n_k}^{-1}, \nu) \to 0, \quad k\to\infty, 
\end{equation*}
whence $\mu^*=\nu$. Hence we have shown that each subsequence 
$\{\mu\circ T_{n_j}^{-1}\}$ contain a further subsequence 
$\{\mu\circ T_{n_{j_m}}^{-1}\}$ that converges weakly to $\nu$. 
Then, by Theorem 2.6 in Billingsley \cite{B}, 
we deduce that $\{\mu\circ T_n^{-1}\}$ converges weakly to $\nu$. 
Denote by $W_2^{2}(\mu_1,\mu_2)$ the 2-Wasserstein distance between $\mu_1$ and 
$\mu_2$. Then we have 
$W_2^{2}(\mu,\nu)=M(T^*)$ and the duality formula 
\begin{equation*}
 W_2(\mu,\nu)=\mathrm{sup}\left\{\int\varphi\,d\mu + \int\psi\, d\nu \right\}, 
\end{equation*}
where the supremum is taken over all bounded continuous functions $\varphi$ and 
$\psi$ such that $\varphi(x) + \psi(y)\le |x-y|^2$, $x,y\in\mathbb{R}^d$. 
See, e.g., Proposition 5.3 in Carmona and Delarue \cite{C}. 
Let $\varepsilon^{\prime\prime}>0$ be arbitrary. Then there exist $\varphi^{\prime}$ 
and $\psi^{\prime}$ such that 
\begin{align*}
 &W_2^{2}(\mu,\nu) \\ 
 &\le \int_{\mathbb{R}^d}\varphi^{\prime}(x)\,d\mu(x) + 
 \int_{\mathbb{R}^d}\psi^{\prime}(y)\,d\nu(y) + \varepsilon^{\prime\prime}. 
\end{align*}
Further, since $\psi^{\prime}$ is bounded and continuous, 
there exists $n_0\in\mathbb{N}$ such that 
$\varphi^{\prime}(x)+\psi(y)\le |x-y|^2$ and 
\begin{equation*}
 \int_{\mathbb{R}^d}\psi^{\prime}(T^*(x))\,d\mu(x) 
 \le \int_{\mathbb{R}^d}\psi^{\prime}(T_n(x))\,d\mu(x) + \varepsilon^{\prime\prime}
\end{equation*}
for $n\ge n_0$. With these choices it follows that for $n\ge n_0$, 
\begin{align*}
 M(T^*)&=W_2^{2}(\mu,\nu) \\ 
 &\le 
 \int_{\mathbb{R}^d}\varphi^{\prime}(x)\,d\mu(x) + 
 \int_{\mathbb{R}^d}\psi^{\prime}(T^*(x))\,d\mu(x) + \varepsilon^{\prime\prime} \\ 
 &\le \int_{\mathbb{R}^d}\left(\varphi^{\prime}(x) + \psi^{\prime}(T_n(x))\right)
  d\mu(x) + 2\varepsilon^{\prime\prime}\notag\\ 
  &\le \int_{\mathbb{R}^d}|x-T_n(x)|^2\,d\mu(x) + 2\varepsilon^{\prime\prime} \\ 
  &\le M_{\lambda_n}(T_n) + 2\varepsilon^{\prime\prime}. 
\end{align*}
Then letting $n\to\infty$ we get 
\begin{equation*}
 M(T^*)\le \liminf_{n\to\infty}M_{\lambda_n}(T_n) + 2\varepsilon^{\prime\prime}.
\end{equation*}
Since $\varepsilon^{\prime\prime}$ is arbitrary, we deduce 
$M(T^*)\le \liminf_{n\to\infty}M_{\lambda_n}(T_n)$. 
On the other hand, \eqref{eq:4} immediately leads to 
$\limsup_{n\to\infty}M_{\lambda_n}(T_n)\le M(T^*)$. 
Therefore, $\lim_{n\to\infty}M_{\lambda_n}(T_n)= M(T^*)$, as wanted. \qed 
\end{Proof}

\section{NUMERICAL EXPERIMENTS}\label{sec:3}
Here we test our two algorithms through several numerical experiments. 
\subsection{Interpolation of synthetic datasets}
All of numerical examples below are implemented in PyTorch on a  Core(TM) i7-13700H with 32GB memory in this subsection.
In these experiments, we describe three experiments on synthetic datasets.
Date size is set to be $500$ for each experiment. 
The Gaussian kernel $K(x,y)=e^{-|x-y|^{2}}$, cost function $c(x,y)=|x-y|^{2}$ 
and the Adam optimizer with learning late $0.0001$ is used.
Here, the function $T(x)$ is described by a multi-layer perception with $1$ hidden layer.
These result that obtained after about $3000$ epochs.
Penalty parameter $\lambda$ defined by $1/\lambda=0.000001$.
\subsubsection{From moon to circle}\label{sec:3.1.1}
\begin{figure}[H]
  \centering
   \includegraphics[width=60mm,scale=0.7]{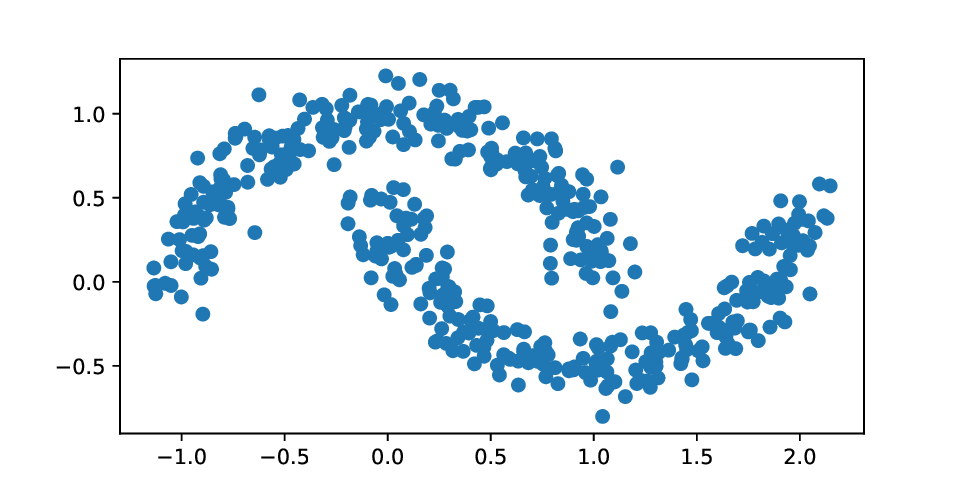}
  \caption{Initial distribution}
  \label{fig:1}
 \end{figure}

\begin{figure}[H]
   \centering
    \includegraphics[width=60mm,scale=0.7]{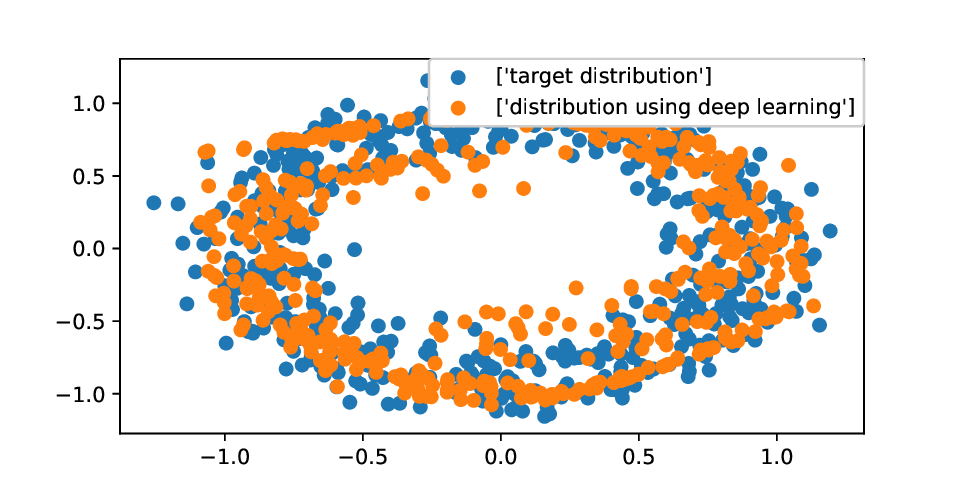}
    \caption{Target distribution (blue) and generated samples (orange). }
    \label{fig:2}
\end{figure}

In this experiment, the initial distribution is 
the well-known synthetic dataset generated by two ``moons" (Figure \ref{fig:1}), 
and the target distribution is the one generated by two ``circles" (Figure \ref{fig:2}). 
\begin{figure}[H]
   \centering
    \includegraphics[width=60mm,scale=0.7]{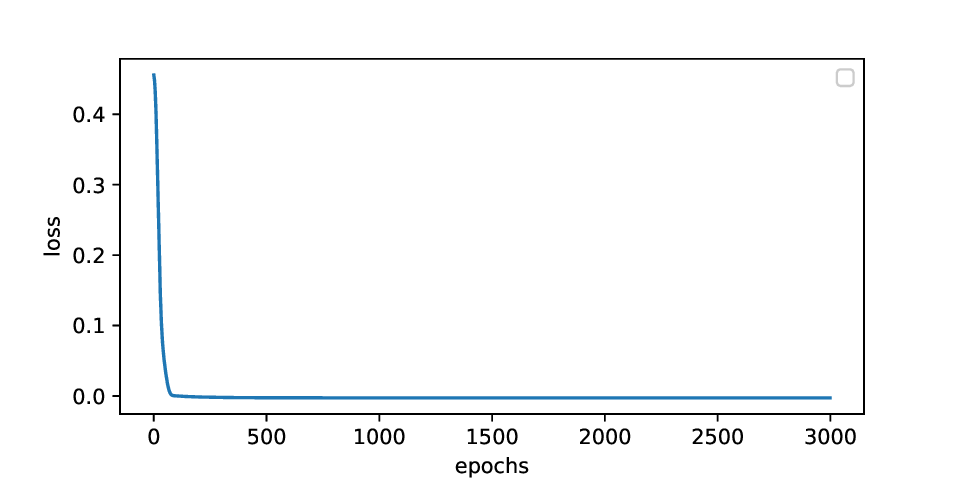}
    \caption{Loss curve}
    \label{fig:3}
\end{figure}
We can see from Figure \ref{fig:2} that the proposed method faithfully generates the target distribution. 
Figure \ref{fig:3} shows the change in loss, and the loss converges in the first 500 epochs.

\subsubsection{From normal distribution to moon}\label{sec:3.1.2}
\begin{figure}[H]
  \centering
   \includegraphics[width=60mm,scale=0.7]{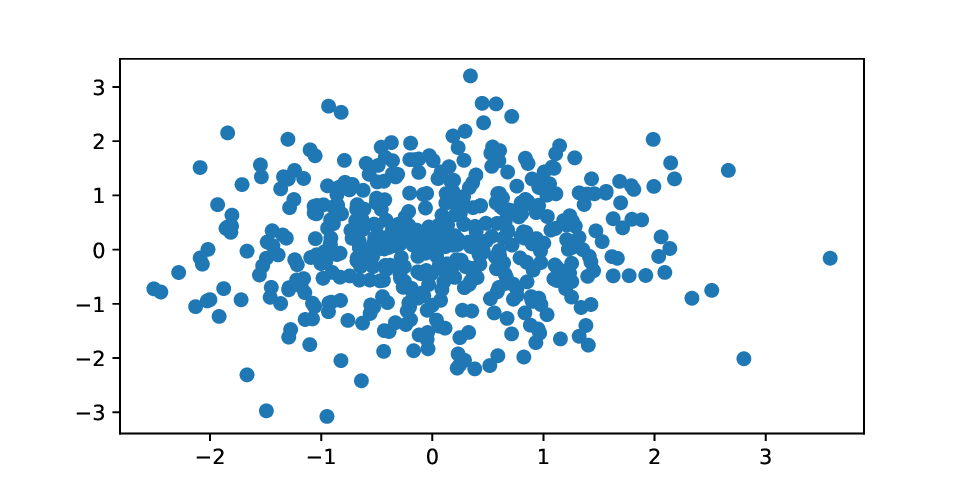}
  \caption{Initial distribution}
  \label{fig:4}
 \end{figure}

\begin{figure}[H]
   \centering
    \includegraphics[width=60mm,scale=0.7]{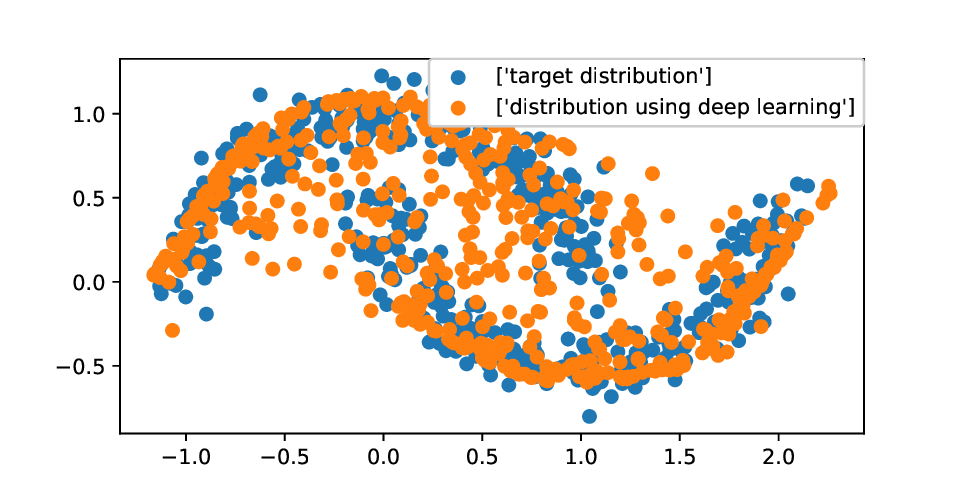}
    \caption{Target distribution (blue) and generated samples (orange)}
    \label{fig:5}
\end{figure}

In this experiment, the initial distribution is a two-dimensional uncorrelated normal distribution with mean 0 and variance 1 (Figure \ref{fig:4}),  
and the final distribution is the synthetic dataset generated by the two moons as in Section \ref{sec:3.1.1} (Figure \ref{fig:5}). 

\begin{figure}[H]
   \centering
    \includegraphics[width=60mm,scale=0.7]{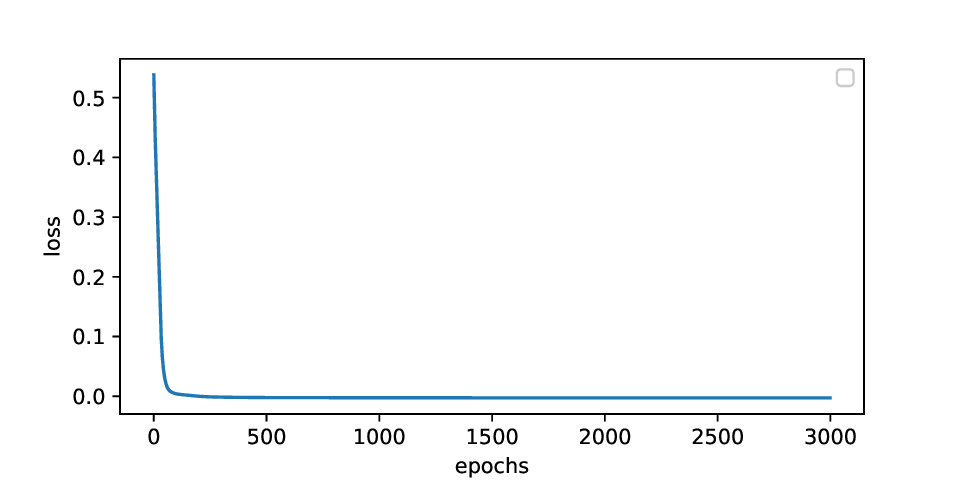}
    \caption{Loss curve}
    \label{fig:6}
\end{figure}
We can see from Figure \ref{fig:5} that the proposed method again generates the target distribution correctly with a small variance. 
Figure \ref{fig:6} shows the change in loss. We can see that the loss converges in the first 500 epochs.

\subsubsection{From normal distribution to normal ditribution}
In this experiment, the initial distribution is a two-dimensional uncorrelated normal distribution with mean $0$ and variance $1$, 
and the target distribution is a two-dimensional uncorrelated normal distribution with mean $5$ and variance $1$.
\begin{figure}[H]
  \centering
   \includegraphics[width=60mm,scale=0.7]{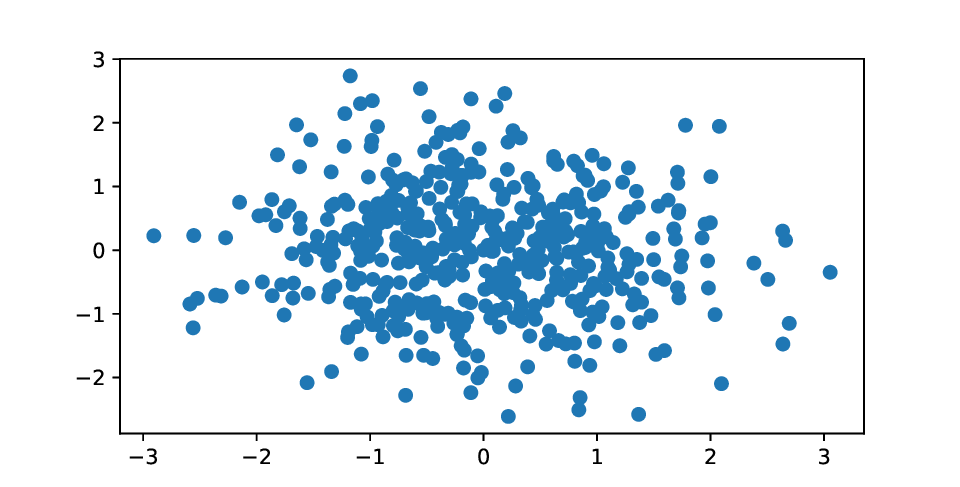}
  \caption{Initial distribution}
  \label{fig:7}
 \end{figure}
\begin{figure}[H]
   \centering
    \includegraphics[width=60mm,scale=0.7]{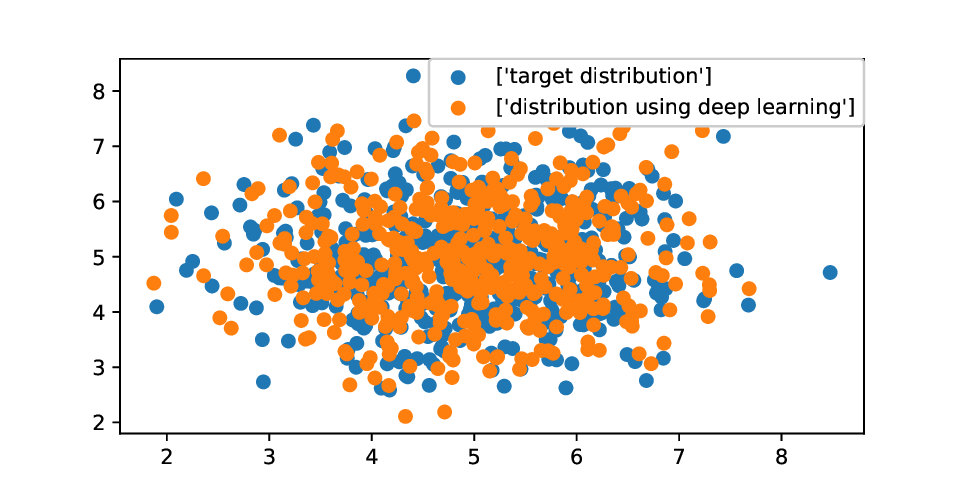}
    \caption{Target distribution (blue) and generated samples (orange)}
    \label{fig:8}
\end{figure}
\begin{figure}[H]
   \centering
    \includegraphics[width=60mm,scale=0.7]{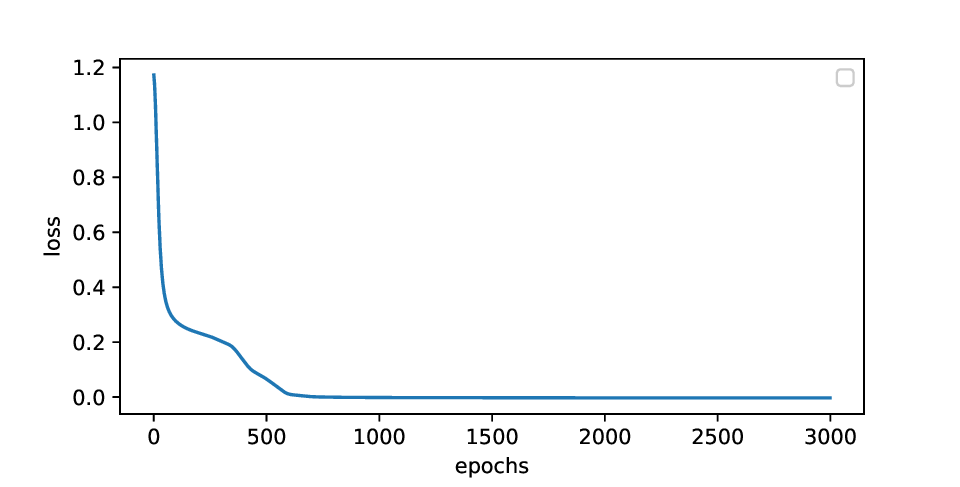}
    \caption{Loss curve}
    \label{fig:9}
\end{figure}
For normal distributions, we see from Figures \ref{fig:7}-\ref{fig:9} that stable generation is achieved as in Sections \ref{sec:3.1.1} and \ref{sec:3.1.2}.  

 \subsection{Comparison with POT}
We compared the performance with POT, the existing Python library as mentioned in Section \ref{sec:introduction}. 
Here, the initial distribution is a two-dimensional uncorrelated normal distribution with mean $0$ and variance $1$, 
and the target distribution is a two-dimensional uncorrelated normal distribution with mean $5$ and variance $1$.
The Gaussian kernel $K(x,y)=e^{-|x-y|^{2}}$, cost function $c(x,y)=|x-y|^{2}$ 
and the Adam optimizer with learning late 0.0001 is used.
Here, the function $T(x)$ is described by a multi\text{-}layer perception with 2 hidden layer.
 Penalty parameter $\lambda$ defined by $1/\lambda=0.000001$.
In this experiment, we compare the number of points that can be calculated, the number of times the calculation can be run, and the accuracy achieved when using an AMD EPYC 9654 with 768GB of memory.
The experimental result is described in Table 1, where SD stands for standard deviation. 
\begin{table}[h]
\centering
\begin{tabular}{ccc} \hline
   Data size & epochs&expectation, SD \\ \hline
   690000&3&4.9707, 1.0154\\ \hline
   350000&5&5.0569, 1.0236\\ \hline
 \end{tabular}
 \caption{The proposed method with CPU}
\end{table}
Here, we set the number of batch size to be $10000$, 
the top row of the table shows the number of training data to be $550000$, 
the number of test data to be $140000$ and so the number of iterations is $55$, 
the bottom row of the table shows the number of training data to be $100000$, 
the number of test data to be $250000$ and so the number of iterations is $10$.
In the upper row, the size of the test data was maintained, 
when the data size was set to be $700000$ (in other words, when the training data was increased), it was not possible to calculate due to CPU usage time.
In the lower row, the size of the training data was maintained, when the size of the test data was set to be $300000$, the calculation became impossible due to the size of the CPU memory.

Next, we perform a similar experiment on NVIDIA H100. 
The experimental result is described in Table 2, where SD again stands for standard deviation. 
\begin{table}[h]
\centering
\begin{tabular}{ccc} \hline
   Data size & epochs &expectation, SD \\ \hline
   610000&85&4.9448, 0.9973\\ \hline
 \end{tabular}
 \caption{The proposed method with GPU}
\end{table}
Here, we set the number of training data to be $600000$, 
the number of test data to be $10000$ and so the number of iterations is $60$. 
The above results show that when the test size was set to be $20000$ and 
the size of the training data was maintained, the calculation became impossible due to the size of the GPU memory.
In addition, when the size of the test data was maintained 
and the number of epochs was set to be $90$, it was not possible to calculate due to GPU usage time.

Next, we compared the calculation speed of the CPU and GPU. 
On both GPU and CPU, we set the number of training data to be $600000$, the number of test data to be $10000$ and so the number of iterations is $60$. 
The epoch number is 1.
\begin{table}[h]
\centering
\begin{tabular}{ccc} \hline
   Processing Unit & time(minutes)&expectation, SD \\ \hline
   CPU&348.4&4.9582, 1.0286\\ \hline
   GPU&17.27&4.8404, 1.0520\\ \hline
 \end{tabular}
 \caption{Comparison of calculation speed on GPU and CPU}
\end{table}
The CPU calculation speed took 20 times longer 
than the GPU calculation speed.

Then, we use the solver {\tt ot.sinkhorn()}\hspace*{0em} in POT  to compare the performance of POT with that of our algorithm on an AMD EPYC 9654 with 768GB of memory. 
The computational complexity of this solver is known to be $O(n^2)$, where $n$ is the input data size.
\begin{table}[h]
\centering
\begin{tabular}{ccc} \hline
   Data size &expectation, SD \\ \hline
   200 &4.8547, 0.9788\\ \hline
   1000 &4.2773, 0.9708\\ \hline
   10000 &3.7472, 0.9393\\ \hline
 \end{tabular}
 \caption{Python Optimal Transport}
\end{table}

 In table 3, calculations were repeated about 10 times for data sizes of $200$, $1000$ and $10000$, and the average values were calculated.

Numerical experiments in this subsection show that our proposed method is a promising option for solving large-scale Monge problems.

\section{\uppercase{CONCLUSION}}\label{sec:4}

In this paper, we derive an approximate solution to the Monge problem by using the embeddings of probability measures and a deep learning algorithm. Through several numerical experiments, we confirmed that our method produces accurate approximate solutions and is efficiently computable on a GPU. In future work, we aim to extend our research to handle larger datasets, explore the Monge problem with alternative cost functions, and investigate numerical solutions for multi-marginal transport problems.

\bibliographystyle{apalike}
\bibliography{smph}

\end{document}